\documentclass{amsart}
\usepackage{amscd, amsmath, amsthm, amssymb}

\newtheorem{theorem}{Theorem}[section]
\newtheorem{lemma}[theorem]{Lemma}
\newtheorem{proposition}[theorem]{Proposition}
\newtheorem{conjecture}[theorem]{Conjecture}

\theoremstyle{definition}     
\newtheorem{definition}[theorem]{Definition}

\newtheorem{example}[theorem]{Example}

\theoremstyle{remark}
\newtheorem{remark}[theorem]{Remark}

\numberwithin{equation}{section}

\newcommand\Supp{\text{\rm Supp}}

\newcommand\Alb{\text{\rm Alb}}
\newcommand\alb{\text{\rm alb}}
\newcommand\Pic{\text{\rm Pic}}  \newcommand\ot{{\otimes}} \newcommand\OO{{\mathcal{O}}}    \newcommand\FF{{\mathcal{F}}}             \newcommand\dimm{{\text{\rm dim}}}  \newcommand\rk{{\text{\rm rank}}}

\begin{document}


\vskip 2pc

\title[A nonvanishing theorem for Q-divisors]
{A nonvanishing theorem for Q-divisors}

\author[Alfred Chen]{Jungkai Alfred Chen}
\address{Department of Mathematics, National Taiwan University,
Taipei, 106, Taiwan}
\address{Mathematics Division, National Center for Theoretical Science at Taipei}
\email{jkchen@math.ntu.edu.tw}

\author[M. Chen]{Meng Chen}
\address{Institute of Mathematics, Fudan University, Shanghai, 200433, P. R. China}
\email{mchen@fudan.edu.cn}

\author[D. -Q. Zhang]{De-Qi Zhang}
\address{Department of Mathematics, National University of Singapore,
2 Science Drive 2, Singapore 117543, Singapore}
\email{matzdq@math.nus.edu.sg}

\thanks {The first author was partially supported by NSC, Taiwan.
The second author was supported by the National Natural Science
Foundation of China (No.10131010). The third author was supported
by an Academic Research Fund of NUS}


\begin{abstract}
We prove a non-vanishing theorem of the cohomology $H^0$ of the adjoint divisor 
$K_X + \lceil L \rceil$ where $\lceil L \rceil$ is the round up of
a nef and big ${\bold Q}$-divisor $L$.
\end{abstract}

\maketitle


\setcounter{section}{0}
\section{ Introduction}

We work over the complex number field {\bf C}. The motivation of
this note is to find an effective version of the famous
non-vanishing theorem of Kawamata and Shokurov (see \cite{KMM}, \cite{Sh}). 
We propose the following:

\begin{conjecture} \label{mainconj}
Let $X$ be a nonsingular projective variety. Let $L$ be a {\bf Q}-divisor on $X$
satisfying the conditions below:
\begin{list}{}{
\setlength{\leftmargin}{10pt}
\setlength{\labelwidth}{6pt}
}
\item[{\rm (1)}] $L$ is nef and big,

\item[{\rm (2)}] $K_X + L$ is nef, and

\item[{\rm (3)}] either $L$ is a Cartier integral divisor, or
$L$ is effective.
\end{list}

Then $H^0(X, K_X + \lceil L \rceil) \ne 0$, where $\lceil L \rceil$ is the round up of $L$.
\end{conjecture}

This kind of non-vanishing problem has been considered by Ambro
[Am], A. Chen-Hacon [CH], Kawamata [Ka], Kollar [Ko], Takayama
[Ta], and others. When $L$ is an integral Cartier divisor,
Kawamata [Ka] has proved the above Conjecture \ref{mainconj} if either $\dim X =
2$, or $\dim X = 3$ and $X$ is minimal (i.e., the canonical
divisor $K_X$ is nef).

Conjecture \ref{mainconj} is slightly different from that of
Kawamata's in \cite{Ka}. It is somewhat general in the sense that the
divisor $L$ in question is not assumed to have integral
coefficients.
It is precisely this non-Cartierness of $L$ that causes a lot of trouble when
estimating $h^0(X, K_X + \lceil L \rceil)$. To elaborate, the
Kawamata-Viehweg vanishing (\cite{KV}, \cite{Vi1})
implies that $h^0(X, K_X + \lceil L \rceil) =
\chi(K_X + \lceil L \rceil)$ when the fractional part of $L$ is of
normal crossings. However, the Riemann-Roch formula for $\chi$ may not be
effective because $\lceil L \rceil$ may not be nef and hence
$\lceil L \rceil . (K_X + \lceil L \rceil)$ may not be
non-negative when $X$ is a nonsingular surface.
The worse thing is that as remarked in a recent
paper of [Xi], there are ${\bf Q}$-Fano 2-folds and 3-folds (see
[Fl]) with vanishing $H^0(X, K_X +(-2K_X))$.

Despite of the observations above, in [Xi] it is proved that
$H^0(X, K_X + (D-K_X)) \ne 0$ for Picard number one Gorenstein del
Pezzo surface $X$ and nef and big ${\bf Q}$-Cartier Weil divisor
$D$. In this note we shall prove the following which is a
consequence of Theorems \ref{theoremk=0}, \ref{ruled},
\ref{theoremk=2} and \ref{theoremk=1} (for the case of integral
Cartier $L$, see [Ka]).

\begin{theorem}\label{main} Let $X$ be a nonsingular projective
surface.
Suppose that $X$ and $L$ satisfy the conditions $(1)$ - $(3)$ in Conjecture \ref{mainconj}.
Suppose further that $X$ is relatively minimal. Then either
$H^0(X, K_X + \lceil L \rceil) \ne 0$ or
$H^0(X, K_X + 4L_{\rm red}) \ne 0$.
\end{theorem}

The second conclusion may occur when $K_X$ is nef
(and the Kodaira dimension $\kappa(X) \ge 1$). In this case,
the conditions in Conjecture \ref{mainconj} are automatically
satisfied whenever $L$ is nef and big. So
if $L$ is an effective {\bf Q}-divisor with all coefficients
less than $1$, then the non-vanishing of 
$H^0(X, K_X + \lceil L \rceil)$ is equivalent to that of $H^0(X,
K_X + L_{\rm red})$, which is stronger than our conclusion. Remark
\ref{remk=1} shows that it is hard to replace the coefficient
$"4"$ in the theorem above by $"1"$.

In Sections 3 and 6 (Theorems \ref{irreg} and Theorem
\ref{ratell}), we prove the following non-vanishing results without assuming
the condition (3) in Conjecture \ref{mainconj},
and the proof presented for the first assertion is applicable to higher dimensional varieties.
The Fourier-Mukai transforms are applied in the proof.

\begin{theorem}
Let $X$ and $L$ be as in Conjecture \ref{mainconj} satisfying the
first two conditions only. Then $H^0(X, K_X + \lceil L \rceil) \ne 0$ if either

\par \noindent
$(i)$ $X$ is a surface with irregularity $q(X) > 0$, or

\par \noindent
$(ii)$ $X$ is a relatively minimal elliptic surface with $\kappa(X) = -\infty$
and $K_X + L$ nef and big.
\end{theorem}

\begin{remark}
$(1)$ In Example \ref{example}, we construct an example of pair
$(X, L)$ satisfying the conditions (1) and (2)
in Conjecture \ref{mainconj} (indeed, both $L$ and $K_X + L$ are nef and big)
but with $H^0(X, K_X + \lceil L \rceil) = 0$.
So an extra condition such as the (3) in Conjecture \ref{mainconj} is
necessary.

\par \noindent
(2) The same example shows that in Kollar's result [Ko] on
non-vanishing of $H^0(X, K_X + M)$ for big divisor $M$, the
``bigness'' assumption on the fundamental group $\pi_1(X)$ is
necessary, because in (1) the $M: = \lceil L \rceil \ge L$ is big
and $\pi_1(X) = (1)$.

\par \noindent
(3) The example also shows the necessity to assume the nefness of
the Cartier integral divisor $D$ (with $(X, B)$ klt and $D - (K_X
+ B)$ nef and big) in Kawamata's conjecture [Ka] for the
non-vanishing of $H^0(X, D)$. Indeed, in the example, we have $\lceil L
\rceil = L + B$ with $B$ a simple normal crossing effective
divisor so that $[B] = 0$, whence $(X, B)$ is klt. To be precise,
let $D : = \lceil L \rceil$. Then $D - (K_X + B) = \lceil L \rceil
- B = L$ is nef and big, $D = K_X + \lceil L \rceil$, and $D$ is
not nef for $D . D_i = -1$ with the notation in the example.
\end{remark}

We end the Introduction with:

\begin{remark} Consider a fibred space $f:V\longrightarrow C$ where
$V$ is a nonsingular projective variety and $C$ a complete curve.
Assume $L$ is a nef and big normal crossing ${\bf Q}$-divisor 
such that $K_V+L$ is nef. The well-known positivity says that
$f_*(\omega_{V/C}\otimes {\mathcal O}_V(\lceil L\rceil))$ is
positive whenever it is not equal to $0$. Pick up a general fibre
$F$ of $f$. The induction of the non-vanishing problem on $F$ may
imply that
$$\text{rk}(f_*(\omega_{V/B}\otimes {\mathcal
O}_V(\lceil L\rceil)))=h^0(F, K_F+\lceil L\rceil|_F)\ge h^0(F,
K_F+\lceil L|_F\rceil)\ne 0.$$

The positivity of $f_*(\omega_{V/B}\otimes {\mathcal O}_V(\lceil
L\rceil))$ has direct applications in studying properties of the
moduli schemes for polarized manifolds. Please refer to \cite{Vi2}
for more details.
\end{remark}
The above remark shows one aspect of the importance of the
effective non-vanishing for {\bf Q}-divisors.

\par \vskip 1pc
{\bf Acknowledgement.} We would like to thank the referee for his / her very
careful reading and suggestions for the improvement of the paper.

\section{Some preparations and an example}\label{prep}

We begin with:

\begin{definition}
A reduced connected divisor $\Gamma$, with only simple normal
crossings, is a rational tree if every component of $\Gamma$ is a
rational curve and the dual graph of $\Gamma$ is a tree (i.e., it
contains no loops).
\end{definition}

Before proving Proposition \ref{treeprop} below, we need two
lemmas in advance.

\begin{lemma}\label{rat}
Let $D = \sum_{j=1}^n D_j$ be a reduced connected divisor
on a nonsingular projective surface $X$. Then $D . (K_X + D) \ge -2$ and the
equality holds if and only if $D$ is a rational tree.
\end{lemma}

\begin{proof}
Note that $\sum_{k < j} D_k . D_j \ge n-1$ and the equality holds
if and only if $D$ is a tree. We calculate:
$D . (K_X+D) = \sum D_j^2 + \sum K_X . D_j + 2 \sum_{k < j} D_k . D_j$
$\ge \sum_j (2p_a(D_j)-2) + 2(n-1) \ge -2$. The lemma follows.
\end{proof}

\begin{lemma}\label{treelem}
Suppose that $X$ is a nonsingular projective surface with
$\chi(\OO_X) = 1$ and $D$ ($\ne 0$) a connected reduced divisor
such that $H^0(X, K_X + D) = 0$. Then the following statements are
true.

\par \noindent
$(1)$ $D$ is a connected rational tree.

\par \noindent
$(2)$ Suppose further that $D$ supports a nef and big divisor
(so $D$ is automatically connected).
Then $\pi_1(X) = (1)$.
\end{lemma}

\begin{proof} 
The Serre duality and Riemann-Roch theorem imply $0 = h^0(X, K_X +
D) = h^1(X, K_X + D) + \frac{1}{2}(K_X + D) . D + \chi(\OO_X) \ge
0 + (-1) + 1$ by Lemma \ref{rat}. Thus $D. (K_X + D) = -2$ and
hence $D$ is a connected rational tree by the same lemma. So
$\pi_1(D) = (1)$. Suppose that $D$ supports on a nef and big
effective divisor. Then the surjective map $\pi_1(D) \rightarrow
\pi_1(X)$ in Nori [No, Cor. 2.3] infers $\pi_1(X) = (1)$.
\end{proof}

The next result is a very important restriction
on $X$ and $L$ in Theorem \ref{main}.

\begin{proposition} \label{treeprop}
Let $X$ be a  nonsingular projective surface with $q(X)=0$ and
$L$ a nef and big effective ${\bold Q}$-divisor such that $H^0(X,
K_X + L_{\rm red}) = 0$. Then $\chi(\OO_X) = 1$, $L_{\rm red}$ is
a connected rational tree and $X$ is simply connected.
\end{proposition}

\begin{proof}
Note that $p_g(X) \le h^0(X, K_X + L_{\rm red}) = 0$.
So $\chi(\OO_X) = 1$.
Now the proposition follows from Lemma \ref{treelem}.
\end{proof}

The result below is used in the subsequent sections.

\begin{lemma}\label{fibk=1}
Suppose that $X$ is a minimal nonsingular projective surface with
Kodaira dimension $\kappa(X) = 1$, $p_g(X) = 0$, and
$\pi_1^{\rm alg}(X) = (1)$
(this is true if $\pi_1(X) = (1)$). Let $\pi : X \rightarrow {\bf
P}^1$ be the unique elliptic fibration with $F$ a general fibre.
The following statements are true:

\par \noindent
$(1)$  $\pi$ has exactly two multiple fibres $F_1, F_2$, and their multiplicities
$m_1, m_2$ are coprime. In particular, if $E$ is horizontal
then $E . F = m_1m_2 m_3$ ($\ge 6$) for some positive integer $m_3$.

\par \noindent
$(2)$ Suppose further that a reduced connected divisor $D$ on $X$ is a rational tree
and contains strictly the support of an effective $\Gamma$ of elliptic fibre type.
Then $\Gamma$ is a full fibre of $\pi$ and of type
$II^*$, $(m_1, m_2) = (2, 3)$ and $E . F = 6$ for some $E$ in $D$
(see \cite{BPV}, Ch V, \S 7, for notation of singular fibres).
\end{lemma}

\begin{proof}
(1) Since $\pi_1(X)^{\rm alg} = (1)$, we have $H_1(X, {\bold Z}) = (0)$
and hence $q(X) = 0$. So $\chi({\OO_X}) = 1$.
Since $\kappa(X) = 1$, there is an elliptic fibration
$\pi : X \rightarrow \pi(X) = {\bf P}^1$, where the image
is ${\bf P}^1$ because $q(X) = 0$.
Let $F_i$ ($1 \le i \le t$) be all multiple fibres of $\pi$,
with multiplicity $m_i$.
If $m = gcd(m_1, m_2) \ge 2$, then the relation
$m(F_1/m - F_2/m) \sim 0$ induces an unramified Galois
${\bf Z}/(m)$-cover of $X$, contradicting the assumption
$\pi_1(X)^{\rm alg} = (1)$.
If $t \ge 3$, then by Fox's solution to Fenchel's conjecture
(see [Fo], [Ch]), there is a base change $B \rightarrow {\bold P}^1$
ramified exactly over $\pi(F_i)$ ($1 \le i \le t$)
and with ramification index $m_i$.
Then the normalization $Y$ of the fibre product $X \times_{{\bf P}^1} B$
is an unramified cover of $X$ (so that the induced fibration $Y \rightarrow B$
has no multiple fibres), again contradicting the assumption that $\pi_1(X)^{\rm alg} = (1)$.

\par
On the other hand, by the canonical divisor formula,
we have $K_X = \pi^*(K_{{\bf P}^1}) + \chi({\OO_X}) F_1
+ \sum_{i=1}^t (m_i - 1) (F_i)_{\rm red}
\sim_{\bf Q} (-1 + \sum_{i=1}^t (1 - \frac{1}{m_i})) F_1$
(so $\pi$ is the only elliptic fibration on $X$).
Since $\kappa(X) = 1$, we see that
$t \ge 2$. Now the lemma follows from the results above.

(2) Since $\Gamma$ is of elliptic fibre type, $0 = K_X . \Gamma = \Gamma^2 = 0$.
Hence $\Gamma$ is a multiple of a fibre of $\pi$.
Since the support of 
$\Gamma$ ($< D$) is a tree, it is of type $I_n^*$, $II^*$, $III^*$ or $IV^*$,
whence $\Gamma$ is a full fibre (and is not a multiple fibre).
By the assumption, there is an $E$ in $D$ such that
$\Supp(E + \Gamma)$ is a connected rational tree.
Thus $E . \Gamma \le 6$ and the equality holds if and only if
$\Gamma$ is of type $II^*$ and $E$ meets the coefficient-6 component of $\Gamma$.
Now (2) follows from (1).
\end{proof}

The example below shows that an assumption like the condition (3)
in Conjecture \ref{mainconj} might be necessary.

\begin{example} \label{example}
We shall construct a nonsingular projective surface $X$ and a {\bf
Q}-divisor $L$ such that the conditions (1) and (2) in Conjecture
\ref{mainconj} are satisfied, but that $H^0(X, K_X + \lceil L
\rceil) = 0$. Indeed, we will see that both $L$ and $K_X + L$ are
nef and big ${\bold Q}$-divisors.

\par
Let $C$ be a sextic plane curve with $9$ ordinary cusps (of type
$(2, 3)$) and no other singularities. This $C$ (regarded as a
curve in the dual plane ${\bf P}^{2*}$) is dual to a smooth plane
cubic (always having $9$ inflectins). Let $\overline{X}
\rightarrow {\bf P}^2$ be the double cover branched at $C$. Then
$\overline{X}$ is a normal $K3$ surface with exactly $9$ Du Val
singularities (lying over the $9$ cusps) of Dynkin type $A_2$. Let
$X$ be the minimal resolution. According to Barth [Ba], these $9
A_2$ are $3$-divisible. That is, for some integral divisor $G$, we
have $3 G \sim \sum_{i=1}^9 (C_i + 2D_i)$ where $\coprod (C_i +
D_i)$ is a disjoint union of the $9$ intersecting ${\bf P}^1$
(i.e., the $9A_2$). Let $H$ be the pull back of a general line
away from the $9$ cusps on $C$. Then $H^2 = 2$ and $H$ is disjoint
from the $9A_2$, so $H . G = 0$. We can also calculate that $G^2 =
-6$. Now let $L = H + G - \frac{1}{3}\sum_{i=1}^9 (C_i + 2D_i)$.
Then $\lceil L \rceil = H + G$ and $\lceil L \rceil^2 = -4$.
Clearly, $K_X + L = L \equiv H$ is nef and big. However, by the
Kawamata-Viehweg vanishing, and Riemann-Roch theorem, we have
$h^0(X, K_X + \lceil L \rceil) = \frac{1}{2} \lceil L \rceil^2 + 2
= 0$.

A similar example can be constructed, if one can find
a quartic surface with $16$ nodes (i.e., a normal Kummer quartic surface).
\end{example}

\section{Irregular surfaces}

In this section, we shall show that Conjecture \ref{mainconj}
holds true (with only the first two conditions there but not the
last condition) for surfaces $X$ with positive irregularity
$q(X)$.

To be precise, let $X$ be a nonsingular projective surface with
$q(X)>0$ and let $\alb: X \to \Alb(X)$ be the Albanese map. Then
we have:

\begin{theorem} \label{irreg}
Let $X$ be a nonsingular projective surface with $q(X)>0$. Let $L$
be a nef and big ${\bf Q}$-divisor such that $K_X+L$ is nef. Then
$H^0(X,K+ \lceil L \rceil) \ne 0$.
\end{theorem}

To see this, we need the following lemma:

\begin{lemma}
Let $\FF \ne 0$ be a $IT^0$ sheaf on an abelian variety $A$, i.e.
for every $i > 0$ we have $H^i(A, \FF \ot P) =0 $ for all $P \in
\Pic^0(A)$. Then $H^0(A,\FF) \ne 0$.
\end{lemma}

The proof can be found in \cite{CH}, but we reprove it here.

\begin{proof}
Suppose on the contrary that $H^0(A, \FF) =0$. Since $\FF$ is
$IT^0$, the Fourier-Mukai transform of $\FF$ is a locally free
sheaf of rank $= h^0(A, \FF)$, hence the zero sheaf. The only
sheaf that transforms to the zero sheaf is the zero sheaf, which
is a contradiction.
\end{proof}

\begin{proof}[Proof of Theorem \ref{irreg}]
Let $f: X' \to X$ be an embedded resolution for $(X,L)$. It is
clear that $f^*L$ is nef and big with simple normal crossing
support. Let $\Delta:=\lceil f^*L \rceil-f^*L$, then $(X', \Delta)$ is
Kawamata log terminal (klt for short; for its definition and 
property, see \cite{KMM}, Def 0-2-10). 
By a property of nef and big  divisor (see e.g. \cite{La}, ex
2.2.17), there is an effective divisor $N$ such that
$A_k:=f^*L-\frac{1}{k} N $ is ample for all $k  \gg 0$. We fix $k$
such that $(X', \Delta+ \frac{1}{k} N)$ is klt. Now we can write $A_k=
(\alb \circ f)^*M +E$ for some ample $\bf{Q}$-divisor $M$ on
$A:=\Alb(X)$ and  effective divisor $E$ on $X'$.  Pick irreducible
divisor $ B \in |(n-1)A|$ for $n \gg 0$  such that $(X', \Delta')$
is klt, where
$\Delta':=\Delta+\frac{1}{k}N + \frac{1}{n}E +\frac{1}{n}B$.
Then we have, where $P' = ({\rm alb} \circ f)^*P$ with
$P \in {\rm Pic}^0(A)$:
$$K_{X'}+\lceil f^*L \rceil + P' \equiv K_{X'}+\frac{(\alb \circ
f)^*M}{n}+\Delta'.$$
 Let $\FF:= \alb_*  f_* \OO_{X'}(K_{X'} +\lceil
f^*L \rceil)$. By Koll\'ar's relative vanishing theorem (cf.
\cite{Ko}, 10.19.2), one sees that $\FF$ is $IT^0$.

We claim that $\FF \ne 0$.

Grant this claim for the time being. By the above lemma, it
follows that $$ h^0(X', K_{X'}+\lceil f^*L \rceil)=h^0(A, \FF )
\ne 0.$$ Since $K_{X'}+\lceil f^*L \rceil= f^*(K_X+\lceil L
\rceil) + \Gamma$, where $\Gamma$ is an exceptional divisor
(possibly non-effective). It's easy to see that $f_*
\OO_{X'}(\Gamma) \subset \OO_X$. By the projection formula, one
has: $$ 0 \ne H^0(X',K_{X'}+\lceil f^*L \rceil)= H^0(X,\OO_X(K_X+\lceil L
\rceil) \ot f_* \OO_{X'}(\Gamma) ) \subset H^0(X,K_X+\lceil L
\rceil).$$ This is the required non-vanishing.

To see the claim, if $\dimm (\alb(X))=2$, then $\alb \circ f$ is
generically finite. Hence it is clear that $\FF \ne 0$. If $\dimm
(\alb(X))=1$. Let $F$ be a general fiber of $\alb \circ f$. Then we have:
$$ \rk(\FF) =h^0(F, (K_{X'} +\lceil f^*L \rceil) |_F ) = h^0(F,
K_F+\lceil f^*L|_F \rceil ).$$ Since $f^*L$ is big, $f^*L.F >0$.
It follows that $\deg(\lceil f^*L|_F \rceil)>0$.

If $g(F) > 0$, then we have $h^0(F, K_F+\lceil f^*L|_F \rceil )
>0$ already. If $g(F)=0$, note that  $K_X+L$ is nef. Note
also that $(K_{X'}+f^*L).F= (K_X+L).f(F)$ since $F$ is general.
This implies that
$$\deg(K_F+\lceil f^*L|_F \rceil) =
(K_{X'}+f^*L).F + (\lceil f^*L \rceil -f^*L).F$$
$$=(K_X+L).f(F) + (\lceil f^*L \rceil -f^*L ).F \ge 0.$$
Hence $h^0(F,K_F+\lceil L_F \rceil )>0$. We conclude that $\FF \ne
0$ and hence the required non-vanishing that $h^0(X, K_X+\lceil L
\rceil) \ne 0$.
\end{proof}

\begin{remark}
In the proof of Theorem \ref{irreg}, without taking log-resolution
at the beginning, one can apply Sakai's lemma \cite{Sa} for
surfaces to get the vanishing of higher cohomology. However, our
argument here works for higher dimensional situation. It shows
that non-vanishing for general fiber gives the non-vanishing.
\end{remark}

\section{Surfaces of Kodaira dimension 0}

In this section, we show that the Conjecture \ref{mainconj} in the Introduction
is true for surfaces $X$ (not necessarily minimal) with Kodaira dimenion $\kappa(X) = 0$.

\begin{theorem}\label{theoremk=0} Suppose that $X$ is a nonsingular projective
surface (not necessarily minimal) of Kodaira dimension $\kappa(X) = 0$. Then
Conjecture \ref{mainconj} is true for effective ${\bf Q}$-divisor $L$.
\end{theorem}

\begin{proof}
By Theorem \ref{irreg}, we may assume that $q(X) = 0$. 
We may also assume that $0 = h^0(X, K_X + L_{\rm red}) $ ($\ge p_g(X)$).
So $X$ is the blow up of an Enriques surface by the classification theory.
On the other hand, $\pi_1(X) = (1)$ by Proposition \ref{treeprop},
a contradiction. This proves the theorem.
\end{proof}

\section{Surfaces with negative $\kappa$, Part I : ruled surfaces}

In this section, we prove Conjecture \ref{mainconj} for relatively minimal surfaces $X$ of Kodaira
dimension $\kappa(X) = - \infty$. By Theorem \ref{irreg}, we may assume that
$q(X) = 0$, so $X$ is a relatively minimal rational surface.
If $X= \bf{P}^2$ or $\bf{P}^1 \times
\bf{P}^1$, it is easy to verify that Conjecture \ref{mainconj} is true
since effective divisor is then nef. We thus assume that $X$ is the Hirzebruch surface
$\mathbb{F}_d$ of degree $d \ge 1$ (though, $\mathbb{F}_1$ is not relatively minimal).

We first fix some notations. Let $\pi:\mathbb{F}_d \to
\bf{P}^1$ be the ruling. Let $F$ be
a general fibre and $C$ the only negative curve
(a cross-section, indeed)
on $\mathbb{F}_d$.  So $C^2 = -d$.

\begin{theorem}\label{ruled}
Let $X$ be a relative minimal surface of Kodaira dimension $\kappa(X) = - \infty$. Then
Conjecture \ref{mainconj} holds for effective ${\bold Q}$-divisor $L$.
\end{theorem}

\begin{proof}
As mentioned above, we assume that $X = \mathbb{F}_d$ for some $d \ge 1$.
Let $L$ be a nef and big effective {\bf Q}-divisor such that $K_X+L$ is nef.
If $\Supp(L)$ does not contain the negative curve $C$, then
$E : = \lceil L \rceil - L$ is effective and nef;
so $\lceil L \rceil = L + E$ is nef and big and
$K_X + \lceil L \rceil = K_X + L + E$ is nef; then
 the Serre duality and Riemann-Roch theorem for Cartier divisor
imply that $h^0(X, K_X + \lceil L \rceil) \ge
\frac{1}{2} \lceil L \rceil (K_X + \lceil L \rceil)
+ \chi(\OO_X) \ge 0 + 1$.
Theorefore, we may assume that $\Supp(L)$ contains $C$.

\par
Write $L = \sum_i c_i C_i + \sum f_j F_j$
where $C_1 = C$, the $C_i$'s are  distinct horizontal components and
$F_j$'s are distinct fibres, where $c_i > 0$, $f_j > 0$.

\par
Suppose on the contrary that $H^0(X, K_X+ \lceil L
\rceil)=0$. Then by Lemma \ref{treelem},
$L_{\rm red}$ is a connected rational tree.
Hence one of the following cases occurs:

\par \vskip 1pc \noindent
Case (i). $L = c_1C_1 + \sum_{j=1}^k f_j F_j$ ($k \ge 0$),

\par \noindent
Case (ii). $L = \sum_{i=1}^k c_i C_i + f_1F_1$ ($k \ge 2$), and
$L_{\rm red}$ is comb-shaped, i.e., $C_i$'s are disjoint
cross-sections.

\par \noindent
Case (iii). $L = \sum_{i=1}^k c_iC_i$ ($k \ge 2$).

\par \vskip 1pc
Recall that $K_X \sim -2 C_1 -(d+2) F$.
The nefness of $K_X + L$ implies:
$$0 \le (K_X+L ). F = -2 + \sum c_i (C_i . F),$$
$$0 \le (K_X + L) . C_1 = d- 2 - d c_1 + \sum_{i \ge 2} c_i (C_i . C_1) + \sum f_j,$$
$$\sum c_i (C_i . F) \ge 2,$$
$$\sum f_j \ge 2 + (c_1-1) d - \sum_{i \ge 2} c_i (C_i . C_1).$$

\par
In Case (i), the above inequalities imply
$c_1 \ge 2$ and $\sum f_j \ge 2 + (c_1 - 1) d \ge d+2$,
whence $\lceil L \rceil = \lceil c_1\rceil C_1 + \sum \lceil f_j \rceil F_j
\sim \lceil c_1 \rceil C_1 + (\sum \lceil f_j \rceil) F
\ge 2 C_1 + (\sum f_j) F \ge 2 C_1 + (d+2) F_1\sim -K_X$.
Hence $H^0(X, K_X + \lceil L \rceil) \ne 0$.

\par
Consider Case (ii). Then one sees easily that
$k = 2$ and $C_2 \sim C_1 + d F_1$ (see [Ha, Chapter V, \S 2]).
By the displayed inequalities, we have
$c_1 \ge 2 - c_2$ and $f_1 \ge 2 + (c_1 - 1) d$.
If $c_2 > 1$ then $\lceil L \rceil \ge C_1 + 2 C_2 + F_1 > -K_X$,
whence $H^0(X, K_X + \lceil L \rceil) \ne 0$.
So we may assume that $c_2 \le 1$.
Then $c_1 \ge 1$ and $f_1 \ge 2$.
Thus $\lceil L \rceil \ge C_1 + C_2 + 2F_1 \sim -K_X$,
whence $H^0(X, K_X + \lceil L \rceil) \ne 0$.

\par
Consider Case (iii). Since $L$ is a connected tree, we may assume
that $C_1 . C_2 = 1$. So $C_2 \sim n(C_1 + d F) + F$ for some
integer $n \ge 1$. Since $C_i$ ($i \ge 2$) is irreducible, we have
$C_i \ge C_1 + d F$ by \cite{Ha}. If $k \ge 3$ or $n \ge 2$, then
we see that $\lceil L \rceil \ge \sum_{i=1}^k C_i >  - K_X$. So
assume that $k = 2$ and $n = 1$. By the inequalities displayed
above, we have $c_1 \ge 2 - c_2$ and $c_2 \ge 2 + (c_1-1) d$. If
$c_2 > 1$ then $\lceil L \rceil \ge C_1 + 2C_2 > -K_X$. So assume
that $c_2 \le 1$. Then $c_1 \ge 2 - 1$ and $c_2 \ge 2 + 0 d$, a
contradiction.
\end{proof}

\section{Surfaces with negative $\kappa$, Part II:  relatively minimal
elliptic}

In this section we consider relatively minimal elliptic
surface $\pi : X \rightarrow B$ with Kodaira dimension $\kappa(X)
= -\infty$. As far as the Conjecture \ref{mainconj} is concerned,
we may assume that the irregularity $q(X) = 0$ by virtue of
Theorem \ref{irreg}. So $X$ is a rational surface and $B = {\bf
P}^1$. By the canonical divisor formula, we see that $\pi$ has at
most one fibre $F_0$ with multiplicity $m \ge 2$; moreover, such
$F_0$ (if exists) is of Kodaira type $I_n$ ($n \ge 0$), and $-K_X
= (F_0)_{\rm red}$.

We show that Conjecture 1.1 is true
if $K_X + L$ is nef and big (but without the assumption of the effectiveness of $L$):

\begin{theorem} \label{ratell}
Let $\pi : X \rightarrow B$ be a relatively minimal elliptic
surface with $\kappa(X) = -\infty$.
Suppose that $L$ is a nef and big ${\bf Q}$-divisor such that $K_X
+ L$ is nef and big. Then $H^0(X, K_X + \lceil L \rceil) \ne 0$.
\end{theorem}

\begin{proof}
By Theorem \ref{irreg}, we may assume that $q(X) = 0$, so
$B = {\bf P}^1$ and $X$ is a rational surface.

\par
Suppose that the ${\bf Q}$-divisor $L$ is nef and big and $K_X + L$ is nef.
Let $F_0 = m (F_0)_{\rm red}$ be the multiple fiber. We set $m = 1$ and let
$F_0$ be a general (smooth) fibre,
if $\pi$ is multiple fibre free. Then $K_X \sim -(F_0)_{\rm red}$.
Let $a > 0$.
Consider the exact sequence:
$$0 \rightarrow {\OO}_X(K_X + \lceil aL \rceil - (F_{\rm 0})_{\rm  red}) \rightarrow
{\OO}_X(K_X + \lceil aL \rceil) \rightarrow {\OO}_{(F_{\rm 0})_{\rm red}}
(K_X + \lceil aL \rceil|_{(F_{\rm 0})_{\rm red}}) \rightarrow 0.$$

Let us find the condition for $aL - (F_{\rm 0})_{\rm red}$ to be nef and big.
Note that $aL - (F_{\rm 0})_{\rm red} \sim aL + K_X = (a-1) L + (K_X + L)$.
So $aL - (F_{\rm 0})_{\rm red}$ is nef and big if either

\par \vskip 1pc
(i) $a > 1$, or

\par
(ii) $a = 1$ and $K_X + L$ is nef and big.

\par \vskip 1pc
Assume that either (i) or (ii) is satisfied.
Then $H^i(X, K_X + \lceil aL \rceil - (F_{\rm 0})_{\rm red}) = 0 =
H^i(X, K_X + \lceil aL \rceil)$ for all $i > 0$,
by Sakai's vanishing for surfaces.
For the integral divisor $M := K_X + \lceil aL \rceil$ and
the reduced divisor $C  := (F_{\rm 0})_{\rm red} $ on $X$,
the above exact sequence implies that
$\chi(\OO_C(M|C) ) = \chi(\OO_X(M)) - \chi(\OO_X(M-C))
= C . M - C . (K_X + C)/2$, where we applied the Riemann-Roch theorem for
both $\OO_X(M)$ and $\OO_X(M - C)$.
Now $C . (K_X + C) = 0$ and $C . M \ge (F_{\rm 0})_{\rm red} . (K_X + aL) > 0$
(for $0 \ne C$ being nef and $K_X + aL$ nef and big), so $\chi(\OO_C(M|C)) > 0$.
By the vanishing above,
$h^0(X, K_X + \lceil aL \rceil) = \chi(\OO_X(M))  = \chi(\OO_X(M-C)) + \chi(\OO_C(M|C) ) =
h^0(X, K_X + \lceil aL \rceil - (F_{\rm 0})_{\rm red}) ) + \chi(\OO_C(M|C) )  > 0  + 0$.
This proves the theorem.
\end{proof}

\begin{remark}
The above argument actually proved the following: let $\pi : X
\rightarrow B$ be a relatively minimal elliptic surface with
$\kappa(X) = -\infty$. Suppose
that $L$ is a nef and big ${\bf Q}$-divisor such that $K_X + L$ is nef.
Then $H^0(X, K_X + \lceil aL \rceil) \ne 0$ provided that either
$a > 1$, or $a = 1$ and $K_X + L$ is nef and big.
\end{remark}

\section{Preparations for surfaces with $\kappa = 1$ or $2$}

Throughout this section, we assume that $X$ is a nonsingular
projective surface with $K_X$ nef and Kodaira dimension $\kappa(X)
= 1$ or $2$. The main result is Proposition \ref{key} to be used in the next section.

\begin{definition}
Up to Lemma \ref{bound}, we let $\Gamma$ be a connected effective integral
divisor on $X$ which consists of smooth rational curves and 
has a (rational) tree as its dual graph.

\par \noindent
$(1)$ We say that $\Gamma$ is of type $A_n'$ (resp. $D_n'$, or $E_n'$)
if its weighted dual graph is of Dynkin type $A_n$ (resp. $D_n$, or $E_n$)
but its weights may not all be $(-2)$.

\par \noindent
$(2)$ 
$\Gamma$ is of type $I_n^*$ (resp. $II^*$, or $III^*$, or $IV^*$)
if $\Gamma$ is of the respective elliptic fibre type
(hence $\Supp(\Gamma)$ is a union of $(-2)$-curves).
$\Gamma$ is of type $I_n^*$' (resp. $II^*$', or $III^*$', or $IV^*$')
if $\Gamma$ is equal to an elliptic fibre of type $I_n^*$ (resp. $II^*$, or $III^*$, or $IV^*$),
including coefficients, but the self intersections of components of $\Gamma$
may not all be $(-2)$.
E.g. $\Gamma = 2\sum_{i=0}^n C_i + \sum_{j=n+1}^{n+4} C_j$ is of type $I_n^*$',
where $C_i + C_0 + C_1 + \dots + C_n + C_j$ is
an ordered linear chain for all $i \in \{n+1, n+2\}$ and $j \in \{n+3, n+4\}$.

\par \noindent
$(3)$ For a divisor $D$ on $X$, we denote by $\#D$ the number of irreducible components
of $D$.

\end{definition}

The assertion(1) below follows from the fact that $C^2 = -2 - C . K_X \le -2$. The others are clear.

\begin{lemma}\label{negdef}
$(1)$ If $C$ is a smooth rational curve on $X$, then $C^2 \le -2$.

\par \noindent
$(2)$ If $\Gamma$ is of type $A_n$', $D_n$' or $E_n$'
then it is negative definite, i.e.,
the intersection matrix of components in $\Gamma$ is negative definite.

\par \noindent
$(3)$ If $\Gamma$ is one of types $I_n^*$, $II^*$, $III^*$ and $IV^*$
(resp. $I_n^*$', $II^*$', $III^*$' and $IV^*$',
but at least one component of $\Gamma$ is not a $(-2)$-curve),
then $\Gamma$ is negative semi-definite (resp. negative definite).

\par \noindent
$(4)$ If $\#\Gamma \le 5$, then $\Gamma$ is negative definite,
unless $\Gamma$ supports a divisor of type $I_0^*$.
\end{lemma}

The Picard number can be estimated in the following way:

\begin{lemma}\label{bound}
Suppose that the
$(-2)$-components of $\Gamma$ do not support a divisor of type $I_0^*$. Let $r
= \min\{9, \, \#\Gamma - 1\}$. Then there is a subgraph $\Gamma'$
 of $r$ components with negative definite intersection matrix.
 In particular, $\rho(X) \ge r+1$. Also if $\rho(X) \le 9$ then $\#\Gamma \le 9$.
  \end{lemma}

\begin{proof}
We have only to prove the
first assertion. By taking a subgraph, we may assume that $\#\Gamma \le 10$.

If $\Gamma$ is a linear chain, then it has negative definite
intersection matrix, and we are done. Thus we may assume that
there exists an irreducible component which meets more than two
other irreducible components.  Let $C_0$ be the irreducible
component that meets $k$ other components with the largest $k$. Then
$\Gamma - C_0$ has exactly $k$ connected components
$\{\Delta_i\}$.
 We may assume that $k \ge 3$. Let
$C_i$ be the irreducible component of $\Delta_i$ that meets $C_0$.

By Lemma \ref{negdef},  if $\#\Delta_i \le 5$ for all $i$ then
each $\Delta_i$ is negative definite. By taking $\Gamma'=\sum
\Delta_i$, we are done.

The remaining cases of $(\#\Delta_1, ..., \#\Delta_k)$ are $\{
(1,1,6), (1,1,7), (1,2,6), (1,1,1,6)\}$. For the case $(1,1,1,6)$,
we take $\Gamma'= \Gamma-C_4$, then now $\Gamma'$ has at least two
connected components: $C_0+C_1+C_2+C_3$ and others. It is clear
that each connected component has at most $5$ irreducible components. Hence $\Gamma'$
is negative definite.  For the cases $(1,1,6)$ and $(1,2,6)$,
similar argument works.

It remains to work with the case $(1,1,7)$. If $C_3$ meets at
least 3 components,
 we take $\Gamma'=\Gamma-C_3$. Then $\Gamma'$ has at least 3
 connected components and each one has length $\le 5$.  If  $C_3$ meets
2 components, say $C_0, C_4$, then we take $\Gamma'= \Gamma-C_4$.
Again, each connected component of $\Gamma'$ has at most $5$ irreducible components .
 This proves the lemma.
\end{proof}

\begin{lemma}\label{elltype}
Suppose that $q(X) = p_g(X) = 0$ and $\pi_1^{\rm alg}(X) = (1)$ 
(these are satisfied in the situation of
Proposition \ref{key}; see its proof).

\par \noindent
$(1)$ We have $\rho(X) \le 10 - K_X^2 \le 10$, and $\rho(X) = 10$ holds
only when $\kappa(X) = 1$.

\par \noindent
$(2)$ For $L$ in Proposition \ref{key}, 
suppose that some $(-2)$-components of $L$ support an
effective divisor $\Gamma$ of elliptic fibre type. Then $\Gamma$
is of type $II^*$, $\kappa(X) = 1$ and $\rho(X) = 10 \le \#L$.
Moreover, $L_{\rm red}$ supports an effective divisor $C$ of type
$I_0^*$' whose central and three of the tip components are all
$(-2)$-curves.
\end{lemma}

\begin{proof}
(1) follows from: $\rho(X) \le b_2(X) = c_2(X) -2 + 4q(X) = 12
\chi(\OO_X) - K_X^2 - 2 = 10 - K_X^2 \le 10$ (Noether's equality).

\par
(2) Since a surface of general type does not contain such
$\Gamma$, we have $\kappa(X) = 1$. By Lemma \ref{fibk=1} and its
notation and noting that $L_{\rm red} > \Supp(\Gamma)$ (for $L$
being nef and big), $\Gamma$ is of type $II^*$ and $\Supp(E +
\Gamma)$ ($\le L_{\rm red}$) supports a $I_0^*$' as described in
(2). Also  $\#L \ge \#\Gamma + 1 = 10$ and $\rho(X) \ge 2 +
(\#\Gamma - 1) = 10$. Thus $\rho(X) = 10$. This proves the lemma.
\end{proof}

By the lemma above and Lemma \ref{bound}, to prove Proposition \ref{key},
we may assume:

\begin{remark}
Assumption: $\#L \le 9$, and the $(-2)$-components of $L$ do not
support a divisor of elliptic fibre type.
\end{remark}

We need three more lemmas in proving Proposition \ref{key}.

\begin{lemma}\label{int}
Let $D = \sum_{i=0}^n D_i$ be a reduced divisor on $X$.
Suppose that $D - D_0$ has a negative definite $n \times n$ intersection matrix
$(D_i . D_j)_{1 \le i, j \le n}$ and $D$ supports a divisor
with positive self intersection.

\par \noindent
$(1)$ We have $\det(D_k . D_{\ell})_{0 \le k, \ell \le n} > 0$
(resp. $< 0$) if $n$ is even (resp. odd).

\par \noindent
$(2)$ Assign formally $G_i : = D_i$ and define
$G_i . G_j := D_i . D_j$ ($i \ne j$)
and $G_i^2 := - x_i$.
Suppose that $(*)$ the $n \times n$ matrix
$(G_i . G_j)_{1 \le i, j \le n}$ is negative definite.
If $G_i^2 \ge D_i^2$ for all $0 \le i \le n$,
then $(**)$ $\det(G_k . G_{\ell})_{0 \le k, \ell \le n} > 0$
(resp. $< 0$) if $n$ is even (resp. odd).

\par \noindent
$(3)$ Suppose that $D_i^2 \le -2$ for all $0 \le i \le n$.
In $(2)$ above for $0 \le k \le n$, choose
the largest positive integer $m_k$ (if exists) such that $(*)$ and $(**)$ in $(2)$ are satisfied
for $G_i$ with $G_k^2 = -m_k$ and $G_i^2 = -2$ ($i \ne k$).
Then $D_k^2 \ge -m_k$.
\end{lemma}

\begin{proof}
For $(1)$, suppose that the matrix in $(1)$ is similar (over {\bf Q}) to
a diagonal matrix $J$. Then the condition implies that $J$
has one positive and $n$ negative diagonal entries. So $(1)$ follows.

For $(2)$, we have only to show that a linear combination of $G_i$ has
positive self intersection.
By the assumption some divisor $\Delta = \sum b_i D_i$ has positive self intersection,
then so is $\Gamma = \sum b_i G_i$ because
$\Gamma^2 = \sum b_i b_j G_i . G_j \ge \sum b_i b_j D_i . D_j
= \Delta^2 > 0$. The (3) follows from (2).
\end{proof}

Let $D = \sum_{i=0}^n D_i$ be a reduced divisor and let
$D = P + N$ be the Zariski decomposition
with $P$ the nef and $N$ the negative part so that
$P$ and $N$ are effective ${\bold Q}$-divisor with $P . N = 0$
(see [Fu1], [Fu2], [Mi]). $D$ supports a nef and big divisor if and only if
$P^2 > 0$.

\par
In Lemmas \ref{zarlem1} and \ref{zarlem2} below, we do not need the bigness of $P$;
in Lemma \ref{zarlem1}, $K_X$ is irrelevant.

\begin{lemma}\label{zarlem1}
$(1)$ Write $P = \sum_{i=0}^n p_i D_i$.
Then $0 \le p_i \le 1$, and $p_i < 1$ holds if and only if $D_i \le \Supp(N)$.

\par \noindent
$(2)$ Write $\Supp(N) = \sum_{i=0}^s D_i$ after relabelling.
Then $(p_0, \dots, p_s)$ is the unique solution of
the linear system
$\sum_{i=0}^n x_i (D_i . D_j) = 0$ ($j = 0, \dots, s$),
where we set $x_j = 1$ ($j > s$).

\par \noindent
$(3)$ Assign formally $G_i := D_i$ and $G_i . G_j = D_i . D_j$ ($i \ne j$).
Suppose that for $\alpha \le i \le \beta$, we assign $G_i^2$ such that
$-2 \ge G_i^2 \ge D_i^2$ and $(G_i . G_j)_{\alpha \le i, j \le \beta}$
is negative definite. Let $(x_i = b_i \,|\, \alpha \le i \le \beta)$ be the unique solution of the
linear system $\sum_{i=0}^n x_i G_i . G_j = 0$ ($\alpha \le j \le \beta$),
where we set $x_j = b_j = p_j$ if $j < \alpha$ or $j > \beta$. Then $b_i \ge p_i$ for
all $\alpha \le i \le \beta$.
\end{lemma}

\begin{proof}
For (1), see [Fu1] or [Mi]. (2) follows from the fact that
$P . D_j = 0$
($0 \le j \le s$) and that $N$ has negative definite (and hence invertible)
intersection matrix.

\par
We prove (3). It suffices to show that
(*) \,\, the sum $\sum_{i=\alpha}^{\beta} (b_i-p_i) G_i . G_j \le 0$
for all $\alpha \le j \le \beta$.

\par
Indeed, write $\sum (b_i - p_i) G_i = A - B$
with $A \ge 0, B \ge 0$ and with no common components in $A$ and $B$;
then the condition (*) implies that
$A . B - B^2 = \sum (b_i - p_i) G_i . B \le 0$;
this and $A . B \ge 0$ and $B^2 \le 0$ imply that $B^2 = 0$
and hence $B = 0$ by the negative-definiteness of $(G_i . G_j)$.

\par
Coming to the sum in (*) above, it is equal to
$\sum_{i=0}^n b_i G_i . G_j - \sum_{i=0}^n p_i G_i . G_j \le 0 - \sum_{i=0}^n p_i D_i . D_j
\le 0$. This proves the lemma.
\end{proof}

\begin{lemma} \label{zarlem2}
Suppose that $\Gamma = D_1+\cdots+D_m$ is an ordered linear chain
contained in $D$ such that $\Gamma . (D - \Gamma) = 1$.
Let $D_t \le \Gamma$ and $D_{m+1} \le D - \Gamma$ such that
$D_t . D_{m+1} = 1$. If either $t = m$ or
$D_t^2 \le -3$, then $\Gamma \le \Supp(N)$.
\end{lemma}

\begin{proof}
Write $P = \sum_j p_j D_j$. If $t = m$,
we set $G_i^2 = -2$ ($1 \le i \le m$) in Lemma \ref{zarlem1}
and obtain $p_i \le b_i = (i/(m+1)) p_{m+1} < 1$ and hence $\Gamma \le \Supp(N)$.
If $D_t^2 \le -3$, we have only to show that $p_t < 1$ because
we already have $p_j < 1$ for every $1 \le j \le m$ with $j \ne t$,
by the previous case. Now
$0 \le P . D_t = p_t D_t^2 + p_{t-1} + p_{t+1} + p_{m+1} < -3 p_t + 3$,
whence $p_t < 1$. This proves the lemma.
\end{proof}

\par
For $L$ in Proposition \ref{key},
let $L_{\rm red} = P + N$ be the Zariski decomposition, so $P \ge
0$ and $N \ge 0$. By the maximality of $P$, we have $L_{\rm red}
\ge P \ge \varepsilon L$ for a suitable small $\varepsilon > 0$
(one can take $\varepsilon$ such that $1/\varepsilon$
is the maximum of coefficients in $L$).
So $P_{\rm red} = L_{\rm red}$. Write
$$P = \sum_{i=0}^n p_i C_i,$$
Then $0 < p_i \le 1$. Note that $p_j = 1$ for some $j$ for
otherwise $\Supp(L) = \Supp(P) \subseteq \Supp(N)$ would be negative
definite. So we assume the following (after relabelling):

\begin{remark} In order to prove the Proposition below,
we may and will assume that $L = P$ and $p_0 = 1$.
\end{remark}

Now we state the main result of the section.

\begin{proposition}\label{key}
Let $X$ be a minimal nonsingular projective surface
(i.e., $K_X$ is nef) with $p_g(X) = 0$.
Suppose that $L$ is a nef and big effective ${\bold Q}$-divisor
supported by a rational tree.
Then $X$ is simply connected and $\Supp(L)$ is connected. Moreover,
either (the number of irreducible components)
$\#L \ge 10 = \rho(X)$ and $\kappa(X) = 1$, or
$\#L \le 9$ and $(A)$ or $(B)$ below is true:

\par \noindent
$(A)$ There is a linear chain $C = \sum_{i=0}^r C_i \le L_{\rm red}$ with
$r \ge 0$ (after relabelling) such that $L_{\rm red}$
$. \sum_{i=0}^r C_i \ge 2$.

\par \noindent
$(B)$ $\Supp(L)$ supports an effective divisor $C$ of type in
$\{$ $I_n^*$', $III^*$', $IV^*$' $\}$ and
the weights of the multiplicity $\ge 2$ components of $C$ are all $(-2)$,
so $C . (K_X + C) = 0$.
Also the type $III^*$' occurs only when $L_{\rm red}$ is given as follows:

\par \noindent
Case $(B1)$. $\kappa(X) = 1$ and $\rho(X) = 10$;
$\det(\Pic(X)) = -1$, and
$\Pic(X)$ is generated by the divisor class of $K_X$ and those of the
$9$ curves in $L_{\rm red} = \sum_{i=0}^8 C_i$;
$C_0$ meets exactly $C_1, C_2, C_3$;
$C_2 + C_4 + C_6$ and $C_3 + C_5 + C_7 + C_8$ are linear chains;
$C_6^2 = -3$ and $C_i^2 = -2$ ($i \ne 6$).
\end{proposition}

\begin{proof}
Since $L$ is nef and big and a rational tree, $\kappa(X) = 1, 2$.
Since $L$ is nef and big, a positive multiple of $L$ is Cartier and 1-connected.
By [No, Cor. 2.3] or the proof of Lemma \ref{treelem},
$\pi_1(X) = (1)$. In particular, $q(X) = 0$ and $\chi(\OO_X) = 1$.

Since $p_0 = 1$ by the additional assumption,
$C_0$ is not in $\Supp(N)$.
Since $0 \le P . C_0 = C_0^2 + \sum p_j$ and $C_0^2 \le -2$,
where $j$ runs in the set so that
$C_j$ meets $C_0$, this $C_0$ meets at least two components of $\Supp(P) - C_0$.
Now the proposition follows from the lemmas below.
\end{proof}

\par \vskip 1pc
By Lemma 7.4, to prove the above Proposition, we 
only need to consider the case $\#P \le 9$.

\begin{lemma} \label{2-comp}
Suppose that $C_0$ meets exactly two components of $\Supp(P) - C_0$.
Then Proposition \ref{key} is true.
\end{lemma}

\begin{proof}
Suppose that $C_0$ meets only $C_1$ and $C_2$ in $\Supp(P) - C_0$.
Then $0 \le P . C_0 = C_0^2 + p_1 + p_2$ and $C_0^2 \le -2$
imply that $p_1 = p_2 = 1$ and $C_0^2 = -2$. Inductively, we can prove that
there is an ordered linear chain (after relabelling) $\sum_{i=a}^b C_i$
in $\Supp(P)$ such that $p_i = 1$ and $C_i^2 = -2$ for all $a \le i \le b$
and $C_a$ (resp. $C_b$) meets $C_{a-1}$ and $C_{a-2}$
(resp. $C_{b+1}$ and $C_{b+2}$) such that $C_{a-2} + C_{a-1} +
2 \sum_{i=a}^b C_i + C_{b+1} + C_{b+2}$ is of type $I_{b-a}^*$'
and Proposition \ref{key} (B) is true.
\end{proof}

\begin{lemma} \label{4-comp}
Suppose that $C_0$ meets at least four components of $\Supp(P) - C_0$.
Then Proposition \ref{key} is true.
\end{lemma}

\begin{proof}
Suppose that $C_0$ meets $C_i$ ($1 \le i \le k$) with $k \ge 4$.
Let $\Delta_i$ ($1 \le i \le k$) be the connected
component of $P_{\rm red} - C_0$ containing $C_i$.
Set $n_i := \#\Delta_i$.
Assume that for only
$1 \le j \le s$ the divisor $C_0 + \Delta_j$ is a linear chain.
By the proof of Lemma \ref{zarlem2}, we have $p_j \le n_j/(n_j+1)$ ($j \le s$).

\par
If $P_{\rm red} . C_0 = C_0^2 + k \ge 2$, then Proposition \ref{key} (A) is true.
So assume that $C_0^2 \le 1-k \le -3$.
Note that $0 \le P . C_0 = C_0^2 + \sum_{i=1}^k p_i
\le C_0^2 + (k-s) + \sum_{i=1}^s p_i \le 1 - s + \sum_{i=1}^s p_i
\le 1 - \sum_{i=1}^s 1/(n_i+1)$.
Suppose that $\#\Delta_i = 1$ for $1 \le i \le s_1$ and $\#\Delta_i \ge 2$
for $i \ge s_1+1$. Then
$$0 \le \sum_{i=s_1+1}^s 1/(n_i+1) \le 1 - s_1/2.$$
Note also that $\#\Delta_j \ge 3$ for all $s+1 \le j \le k$.
Thus,
$$3k-s-s_1 = s_1+2(s-s_1)+3(k-s) \le \#P - 1 \le 8.$$
These two highlighted inequalities imply that
$s = 2$ and $(\#\Delta_1,\dots,\#\Delta_k) = (1, 1, 3, 3)$.

\par
Note that $C_0$ meets the mid-component $C_j$ of $\Delta_j$ ($j = 3, 4$).
By the proof of Lemma \ref{zarlem2}, for every $j$ with $j \ne 0, 3, 4$,
we have $p_j \le 1/2$.
Thus $0 \le P . C_0 =
C_0^2 + \sum_{i=1}^4 p_i \le -3 + (1/2)+(1/2)+1+1 = 0$,
so $C_0^2 = -3$ and $p_3 = p_4 = 1$.
Now $0 \le P . C_3 \le C_3^2 + p_0 + (1/2) + (1/2)$
implies $C_3^2 = -2$.
So $P_{\rm red} . (C_0+C_3) = 2$ and Proposition \ref{key} (A) is true.
\end{proof}

Now we assume that $C_0$ meets exactly three components $C_i$ ($i = 1, 2, 3$)
of $\Supp(P) - C_0$. Let $\Delta_i$ be the connected component of $\Supp(P) - C_0$
containing $C_i$. Set $n_i := \#\Delta_i$.
Then $\sum_{i=1}^3 n_i = \#P - 1 \le 8$.
We may assume that $n_1 \le n_2 \le n_3$.
Then $n_3 \le 6$ and $n_1 \le 2$, so $C_0 + \Delta_1$
is a linear chain. By the proof of Lemma \ref{zarlem2}, we have
$p_1 \le n_1/(n_1+1) < 1$. This and
$0 \le P . C_0 = C_0^2 + p_1 + p_2 + p_3$, together with $C_0^2 \le -2$,
imply that $C_0^2 = -2$. We shall apply Lemma \ref{int} frequently,
where $G_0$ can be chosen as $C_0$ or $C_3$.

\begin{lemma}\label{1-1}
Suppose that $\#\Delta_i = 1$ for $i = 1$ and $2$
(this is true if $\#\Delta_3 = 6$). Then Proposition \ref{key} is true.
\end{lemma}

\begin{proof}
By the proof of Lemma \ref{zarlem2}, we have $p_i \le 1/2$ for $i = 1$ and $2$.
Now $0 \le P . C_0 = C_0^2 + p_1 + p_2 + p_3$ (and $C_0^2 = -2$) imply
$p_3 = 1$ and $p_i = 1/2$ ($i = 1, 2$).
By Lemma \ref{2-comp} and \ref{4-comp} (applied to $C_3$),
we may assume that $C_3$ meets exactly three components $C_0, C_4, C_5$ of
$\Supp(P) - C_3$. If $C_3^2 = -2$, then $P_{\rm red} . (C_0 + C_3) = 2$
and Proposition \ref{key} (A) is true.
Suppose that $C_3^2 \le -3$. Then as above $C_3^2 + p_0 + p_4 + p_5 = P . C_3 \ge 0$
implies that $C_3^2 = -3$ and $p_4 = p_5 = 1$.
(Of course, $p_0 = 1$ is always assumed). Again by the same Lemmas
we may assume that $C_i$ ($i = 4, 5$) meets exactly three components
(one of which is $C_3$). Then $\#P \ge 10$, a contradiction to the
additional assumption $\# P \le 9$.
\end{proof}

\begin{lemma}\label{linchain}
Suppose that $C_0 + \Delta_i$ is a linear chain for all $i = 1, 2, 3$.
Then Proposition \ref{key} is true.
\end{lemma}

\begin{proof}
Note that  $\sum_{i=1}^3 n_i = \#P - 1 \le 8$.
By Lemma \ref{1-1}, we may assume that $n_3 \le 5$.
Except the cases below, $P$ is negative definite or semi-definite by Lemma \ref{negdef},
which is impossible:
$$(n_1, n_2, n_3) = (1, 3, 4), (2, 2, 4), (2, 3, 3), (2, 2, 3).$$

In the first (resp. the last three) cases, $\Supp(P)$ supports a divisor $D$ of type $III^*$'
(resp. $IV^*$'). We need to show that the coefficient $\ge 2$ components
of $D$ are $(-2)$-curves and that $P_{\rm red} = L_{\rm red}$
is given as in Proposition \ref{key} (B1) in the first case.
These follow from Lemma \ref{int} applied to all $0 \le k \le 8$.
For instance, in notation of Proposition \ref{key} (B1),
if we set $-2 \ge G_k^2 = - x_k$ ($k = 6, 8$) and $G_j^2 = -2$ ($j \ne 6, 8$),
then $\det(G_i . G_j)_{0 \le i, j \le 8} = -4 + 3x_6 + 4x_8 - 2x_6x_8 > 0$
provided that $G_k^2 \ge C_k^2$; also if the case $C_i^2 = -2$ ($i \ne 8$)
occurs then the $(-2)$-components of $P$ support a $C$ of {\it elliptic fibre type} $III^*$.

When $P_{\rm red}$ is as in Proposition \ref{key} (B1),
one can check that the lattice ${\bold Z}[K_X, C_i's]$ generated by
the divisor class of $K_X$ and those of the nine curves in $P$,
has determinant $K_X^2 - 1$. Note also that $\rho(X) \le 10 - K_X^2$.
So either $K_X^2 = 0$ (and $\kappa(X) = 1$),
$\Pic(X) = {\bold Z}[K_X, C_i's]$ (noting that $\Pic(X)$ is torsion free
for $\pi_1(X) = (1)$), $\det(\Pic(X)) = -1$ and
Proposition \ref{key} (B1) is true, or $K_X^2 = 1$;
but the latter situation implies, after a direct calculation, that
$K_X$ is numerically (and hence linearly, for $\pi_1(X) = (1)$)
equivalent to an effective integral divisor $\sum k_i C_i$
with $(k_0, k_1, \dots, k_8) = (10, 5, 7, 8, 4, 6, 1, 4, 2)$,
contradicting the assumption that $p_g(X) = 0$.
\end{proof}

\begin{lemma}\label{the5}
Suppose that $n_3 = \#\Delta_3 = 5$.
Then Proposition \ref{key} is true.
\end{lemma}

\begin{proof}
Since $n_1 + n_2 = \#P -1 - n_3\le 3$, we have $(n_1, n_2) = (1, 1), (1, 2)$
and $C_0 + \Delta_i$ ($i = 1, 2$) is a linear chain.
By Lemma \ref{linchain}, we may assume that
$C_0 + \Delta_3$ is not a linear chain.

\par
We shall apply Lemma \ref{int} to deduce the result.
The case $\#P \le 8$ can be reduced to the case $\#P = 9$ because
if an effective $P_1$ with $\#P_1 = 8$ supports a nef and big divisor
then $P$ with $P > P_1$ supports a nef and big divisor too.
So $(n_1, n_2) = (1, 2)$.

\par
Suppose that $\Delta_3$ is a linear chain.
By Lemma \ref{int}, we have $C_3^2 = -2$,
whence $P_{\rm red} . (C_0 + C_3) = 2$ and Proposition \ref{key} (A) is true.
Indeed, if we set $-2 \ge G_3^2 = -x_3$ and $G_j^2 = -2$ ($j \ne 3$)
then $0 < \det(G_{i, j})_{0 \le i, j \le 8}$ equals $114 - 45 x_3$
(when $C_0$ meets the middle component of $\Delta_3$), or
$98 - 40 x_3$ (otherwise), provided that $G_3^2 \ge C_3^2$
(to guarantee the inequality on the determinant).

\par
Suppose that $\Delta_3$ is not a linear chain. Then it is of type
$I_0^*$' or $D_5$'. We denote by $C_{\ell}$ the central component.
Consider the case where $\Delta_3$ is of type $I_0^*$'.
If $C_3$ in $\Delta_3$ (and meeting $C_0$) is a tip component (resp. the central component $C_{\ell}$)
of $\Delta_3$, then applying Lemma \ref{int}, we have $C_{\ell}^2 = -2$
(resp. $C_{\ell}^2 = -2, -3$).
Thus $P_{\rm red} . C_{\ell} \ge 2$ and Proposition \ref{key} (A) is true.

Consider the case where $\Delta_3$ is of type $D_5'$
so that $C_{\alpha} + C_{\beta} + C_{\ell}$ is
the ordered linear chain in $\Delta_3$.
If $C_3$ is $C_{\alpha}$ (resp. $C_{\beta}$, or $C_{\ell}$, or a tip component
$C_{\gamma} \ne C_{\alpha}$ of $\Delta_3$),
applying Lemma \ref{int}, we have $C_i^2 = -2$ for all $C_i$ in $C$
so that $P_{\rm red} . C = 2$ and hence Proposition \ref{key} (A) is true, where $C$ equals
$C_0 + C_{\alpha} + C_{\beta} + C_{\ell}$ (resp. $C_0 + C_{\beta}$, or $C_{\ell}$,
or $C_0 + C_{\gamma} + C_{\ell}$).
\end{proof}

\begin{lemma}\label{the4}
Suppose that $n_3 = \#\Delta_3 = 4$.
Then Proposition \ref{key} is true.
\end{lemma}

\begin{proof}
As in the previous lemma, we only need to consider the case $\#P = 9$.
Then $n_1 + n_2 = \#P - 1 - n_3 = 4$ and $(n_1, n_2) = (1, 3), (2, 2)$.
So $C_0 + \Delta_1$ is a linear chain.

\par
Consider the case that $C_0 + \Delta_2$ is not a linear chain.
Then $(n_1, n_2) = (1, 3)$. If $C_2^2 = -2$, then
$P_{\rm red} . (C_0 + C_2) = 2$ and Proposition \ref{key} (A) is true.
Suppose that $C_2^2 \le -3$. By Lemma \ref{zarlem2}, $\Delta_1 + \Delta_2
\le \Supp(N)$, and by Lemma \ref{zarlem1} with
$G_1^2 = -2$ (resp. $G_2^2 = -3$)
we have $p_1 \le 1/2$ (resp. $p_2 \le 1/2$ , and the other two components of $\Delta_2$
have coefficients less than or equal to $1/4$ in $P$).
This and $0 \le P . C_0 = C_0^2 + p_1 + p_2 + p_3$ imply that
$p_3 = 1$. By Lemma \ref{2-comp} and \ref{4-comp} we may assume that
$C_3$ meets exactly three components (one of which is $C_0$),
so $\Delta_3$ is a linear chain.
If $C_3^2 = -2$, then $P_{\rm red} . (C_0 + C_3) = 2$
and Proposition \ref{key} (A) is true.
If $C_3^2 \le -3$, then $\Delta_3 \le \Supp(N)$ by Lemma \ref{zarlem2};
applying Lemma \ref{zarlem1} with $G_3^2 = -3$, we have
$p_3 \le 6/11$ (the coefficients of components of $\Delta_3$ in
$P$ are respectively less than or equal to $2/11, 4/11, 6/11, 3/11$);
this leads to that $P . C_0 < 0$, a contradiction.

\par
Therefore, we may assume that $C_0 + \Delta_2$ is a linear chain
but $C_0 + \Delta_3$ is not a linear chain (see Lemma \ref{linchain}).
If $\Delta_3$ is a linear chain and
$C_3^2 = -2$, then $P_{\rm red} . (C_0 + C_3) = 2$
and Proposition \ref{key} (A) is true.
If $\Delta_3$ is a linear chain and $C_3^2 \le -3$,
then as above we have $p_3 \le 6/11$
and $p_1 + p_2 \le \sum_{i=1}^2 \frac{n_i}{n_i+1} \le 4/3$.
This leads to that $0 \le P . C_0 = C_0^2 + p_1 + p_2 + p_3
\le -2 + (4/3) + (6/11) < 0$, a contradiction.

\par
Thus we may assume that $\Delta_3$ is not a linear chain, hence
of type $D_4$' with the central component $C_{\ell}$.
For both cases of $(n_1, n_2) = (1, 3)$ and $(2, 2)$,
if  $C_3$ is a tip component
(resp. $C_{\ell}$) of $\Delta_3$, then applying Lemma \ref{int}
we have $C_i^2 = -2$ for all $C_i$ in $C$ so that $P_{\rm red} . C = 2$
and Proposition \ref{key} (A) is true, where
$C$ equals $C_0 + C_3 + C_{\ell}$ (resp. $C_{\ell}$).
This proves the lemma.
\end{proof}

\begin{lemma}\label{the3}
Suppose that $n_3 = \#\Delta_3 \le 3$.
Then Proposition \ref{key} is true.
\end{lemma}

\begin{proof}
As in the previous lemmas, we may assume that $\#P = 9$,
so $(n_1, n_2, n_3) = (2, 3, 3)$.
Hence $C_0 + \Delta_1$ is a linear chain.
By Lemma \ref{linchain}, we may assume that $C_0 + \Delta_3$ is not a linear chain.
When $C_0 + \Delta_i$ ($i = 2$ or $3$) is not a linear chain and $C_i^2 = -2$,
we have $P_{\rm red} . (C_0 + C_i) = 2$
and Proposition \ref{key} (A) is true.
So assume that $C_3^2 \le -3$.
Then $p_3 \le 1/2$ by Lemmas \ref{zarlem2} and \ref{zarlem1} with $G_3^2 := -3$.
Also we may assume either $C_0 + \Delta_2$
is a linear chain or otherwise and $C_2^2 \le -3$ (and hence $p_2 \le 1/2$).
If the former case occurs, by the proof of Lemma \ref{zarlem2}, we
have $p_i \le n_i/(n_i+1)$ ($i = 1, 2$) and $0 \le P . C_0 =
C_0^2 + p_1 + p_2 + p_3 \le -2 + (2/3) + (3/4)  + (1/2) < 0$, a contradiction.
If the latter case occurs, then $0 \le P . C_0 \le -2 + (2/3) + (1/2) + (1/2) < 0$,
a contradiction. This proves the lemma.
The proof of Proposition \ref{key} is also completed.
\end{proof}

\section{Surfaces of Kodaira dimension 1 or 2}

In this section we shall prove the two theorems below:

\begin{theorem} \label{theoremk=2}
Let $X$ be a minimal nonsingular projective surface of Kodaira dimension $2$.
Let $L$ be a nef and big effective ${\bf Q}$-divisor.
Then $H^0(X, K_X + 3 L_{\rm red}) \ne 0$.
\end{theorem}

\begin{theorem} \label{theoremk=1}
Let $X$ be a minimal nonsingular projective surface of Kodaira dimension $1$.
Let $L$ be a nef and big effective ${\bf Q}$-divisor.

\par \noindent
$(1)$ We have $H^0(X, K_X + 4 L_{\rm red}) \ne 0$.

\par \noindent
$(2)$ Suppose that $H^0(X, K_X + 3 L_{\rm red}) = 0$.
Then $L_{\rm red}$ contains at least its name sake with
$9$ components given in Proposition \ref{key} $(B1)$.
Further, $\pi_1(X) = (1)$, $\rho(X) = 10$, $\det(\Pic(X)) = -1$ and
the elliptic fibration $\pi : X \rightarrow {\bf P}^1$ has
exactly two multiple fibres, and their multiplicities are $2$ and $3$.
The $\Pic(X)$ is generated by the divisor class of $K_X$ and those of the $9$
components of $L$.
\end{theorem}

\par
We now prove Theorems \ref{theoremk=2} and \ref{theoremk=1} simultaneously.
By Theorem \ref{irreg}, we may assume that $q(X) = 0$.
We may also assume that $H^0(X, K_X + L_{\rm red}) = 0$,
so $p_g(X) = 0$ and $\chi({\OO}_X) = 1$.
By Proposition \ref{treeprop},
the $L_{\rm red}$ is a connected rational
tree and $\pi_1(X) = (1)$. So we can apply Proposition \ref{key}.

\par \vskip 1pc
Consider first the case $\#L \le  9$
(this is true if $\kappa(X) = 2$ by Proposition \ref{key}).
We apply Proposition \ref{key}.
If Proposition \ref{key} (A) occurs, applying the Serre duality and Riemann Roch theorem,
we have
$h^0(X, K_X + L_{\rm red} + C) \ge \frac{1}{2}(K_X + L_{\rm red} + C) . (L_{\rm red} + C)
+ \chi(\OO_X) =
\frac{1}{2}\{(K_X + L_{\rm red}) . L_{\rm red} + (C^2 + K_X . C) + 2 C . L_{\rm red}\}
+ 1 \ge \frac{1}{2}\{(-2) + (-2) + 2 \times 2\} + 1 = 1$,
where the terms $(-2)$ are due to the fact
that both $L_{\rm red}$ and $C$ are connected rational trees.
Since $2 L_{\rm red} \ge L_{\rm red} + C$, the theorems follow in this case.

\par
Suppose Proposition \ref{key} (B) occurs.
As above we have $h^0(X, K_X + C) \ge \frac{1}{2}(K_X + C) . C
+ \chi(\OO_X) = 0 + 1 = 1$.
If $C$ is of type $III^*$', then $L_{\rm red}$ is given in
Proposition \ref{key} (B1) (so $\kappa(X) = 1$) and we have $4 L_{\rm red} \ge C$;
thus both Theorems \ref{theoremk=1} and \ref{theoremk=2} are true
by Lemma \ref{exceptional} below.
If $C$ is of other type, then $3 L_{\rm red} \ge C$.
This proves the theorems.

\par \vskip 1pc
It remains to consider the case where $\#L \ge 10$. So $\kappa(X) = 1$
and $\rho(X) = 10$ by Proposition \ref{key}.
By Lemma \ref{elltype} and the calculation above, we may
proceed with the additional assumption that no divisor of
elliptic fibre type is supported by
some $(-2)$-components of $\Supp(L)$.
By Lemma \ref{bound}, we have $\rho(X) = 10$
and we may assume that $\Pic(X) \otimes {\bold Q}$ is generated by
$C_i$ ($1 \le i \le 10$) in $L_{\rm red}$ after relabelling:
first find 9 components of $L_{\rm red}$ having a negative definite intersection matrix,
and then the $10th$ generator can be found from $\Supp(L)$ because $L$ is
nef and big (so not negative definite).

\par
Theorefore, $K_X$ is numerically equivalent to a ${\bf Q}$-linear combination of
$C_i$ ($1 \le i \le 10$). Split the combination as $L_2 - L_1$ so that
$K_X + L_1 \sim_{\bf Q} L_2$, where both $L_j$ are effective, $(L_{j})_{\rm{red}} \le \sum_{i=1}^{10} C_i$
and there is no common component of $L_1$ and $L_2$.
Since $\kappa(X) = 1$, we have $L_2 > 0$. Also $L_2$ is nef, noting that $K_X$ is nef.

\par
Suppose that $L_2$ is not big. Then $0 = L_2^2 = L_2 . K_X + L_2 . L_1 \ge 0 + 0$.
Thus $K_X . L_2 = 0$ and hence $L_2$ is contained in fibres of the
elliptic fibration $\pi : X \rightarrow {\bold P}^1$, noting that $q(X) = 0$,
(so that $K_X$ is numerically equal to a positive multiple of a fibre).
This and the fact that $L_2^2 = 0$ and fibre components are negative semi-definite
[Re],
imply that $L_2 = \sum b_j F_j$ where $b_j$'s are positive rational numbers
and $F_j$'s are full fibres, whence $(-2)$-components of
$L_{\rm red}$ ($\ge (L_{2})_ {\rm{red}}$) support an elliptic fibre, contradicting
the additional assumption.

\par
Therefore, $L_2$ is nef and big.
Thus $L_1 \ne 0$ because $K_X \sim_{\bf Q} L_2 - L_1$ is nef but not big.
This and the fact that $\#L_1 + \#L_2 = \#(L_1 + L_2) \le 10$ imply
that $\#L_2 \le 9$.

\par
So we are reduced to the case $\#L \le 9$ after replacing $L_{\rm red}$
by its subdivisor $(L_{2})_{\rm{red}}$.
This proves the theorem.

\begin{lemma}\label{exceptional}
Suppose that $X$ is a minimal nonsingular projective surface
of Kodaira dimension $\kappa(X) = 1$ and $p_g(X) = 0$. Let $D$ be
the reduced divisor given in Proposition \ref{key} $(B1)$ (denoted as $L_{\rm red}$ there).
Then the elliptic fibration $\pi : X \rightarrow {\bf P}^1$ has
exactly two multiple fibres, and their multiplicities are $2$ and $3$.
\end{lemma}

\begin{proof}
We change the labelling and write $D = \sum_{i=0}^8 D_i$,
where $D_0$ meets $D_1, D_5$ and $D_8$;
$D_1 + \dots + D_4$ amd $D_5 + D_6 + D_7$ are linear chains;
$D_7^2 = -3$ and $D_i^2 = -2$ ($i \ne 7$).
We can check that $D$ supports a nef and big divisor
(the Zariski positive part of $D$):
$P = D_0 + \frac{1}{5}(D_4 + 2D_3 + 3D_2 + 4D_1) +
\frac{1}{7}(D_7 + 3D_6 + 5D_5) + \frac{1}{2} D_8$.
Indeed, $P^2 = P . D_0 = 1/70$.
By [No, Cor. 2.3] or the proof of Lemma \ref{treelem}, we have $\pi_1(X) = (1)$,
whence $q(X) = 0$ and $\chi(\OO_X) = 1$.
As in the proof of Lemma \ref{elltype}, we have
$\rho(X) \le 10$.
We can check that the lattice ${\bf Z}[K_X, D_i's]$
generated by the divisor classes of $K_X$ and those of the $9$ curves of $D$
has determinant $-1$. So this lattice equals $\Pic(X)$ and $\rho(X) = 10$,
noting that $\Pic(X)$ is torsion free for $\pi_1(X) = (1)$.

\par
By Lemma \ref{fibk=1} (and the notation there)
and by the canonical divisor formula
we have $K_X \sim_{\bf Q} (1 - \frac{1}{m_1} - \frac{1}{m_2}) F_1$.
We still have to show that $(m_1, m_2) = (2, 3)$.
Let $F_3$ be the fibre of $\pi$ containing the eight $(-2)$-components of $D$.
Then $F_3$ must be of type $II^*$, so there is a
$(-2)$-curve $G$ such that $G$ and the eight $(-2)$-components of $D$
support the fibre $F_3$ (whence $G . D_4 = 1$ and $G . D_i = 0$ ($i \ne 4, 7$)).

\par
On the other hand,
express $G \sim k K_X + \sum_{i=0}^8 d_iD_i$ for some
integers $k, d_i$. Intersecting the equality by $K_X$,
we obtain $0 = d_7 D_7 . K_X = d_7$.
So $k K_X \sim G - \sum_{i \ne 7} d_i D_i$ and the RHS
is supported on the fibre $F_3$ and has self intersection $0$ (because $K_X^2 = 0$).
Since the fibre components are negative semi-definite,
this implies that the RHS is a multiple of $F_3$.
Now $G$ has coefficient $1$ in $F_3$, so the RHS $= F_3$.
Namely, $kK_X \sim F_3$, or $K_X \sim_{\bf Q} F_3/k$.
Comparing with the expression of $K_X$ in the previous paragraph,
we obtain: $\frac{1}{k} = (1 - \frac{1}{m_1} - \frac{1}{m_2})$.
Simplifying, we obtain: $m_1 m_2 = k(m_1 m_2 - m_1 - m_2)$.
Since $m_1$ and $m_2$ are coprime, we have $m_1 m_2 | k$.
So $k = m_1 m_2$ and $m_1 m_2 - m_1 - m_2 = 1$,
or $1 = \frac{1}{m_1} + \frac{1}{m_2} + \frac{1}{m_1m_2}$.
One sees then $(m_1, m_2) = (2, 3)$.
By the way, then $F_3 \sim F_1 \sim 6 K_X$.
Intersecting this relation with $D_7$, we see that $D_7$
is a $6$-section and $D_7 . G = 4$. This proves the lemma.
\end{proof}

\begin{remark} \label{remk=1}
The non-vanishing of $H^0(X, K_X + L_{\rm red})$ or
$H^0(X,$ $K_X + \lceil L \rceil)$, when $\kappa(X) = 1$,
is subtle and is not easy to be proven at all.
Indeed, suppose that $X$ is a minimal nonsingular projective surface
with Kodaira dimension $1$, $q(X) = 0$ and $p_g(X) = 0$.
Let $\pi : X \rightarrow {\bf P}^1$
be the elliptic fibration.
Suppose that there is a
type $II^*$ elliptic fibre $F_0$ and also there is a $6$-section $E$ ($\cong {\bf P}^1$)
such that $E$ meets the multiplicity-6 component of $F_0$.
(We have this possible situation in mind:
$\pi$ has exactly two multiple fibres. Their multiplicities
are 2, 3; see Lemma \ref{fibk=1}).
Then $L = \frac{1}{6n}(E + nF_0)$ is nef and big for $n > > 0$.
Clearly, $L_{\rm red}$ is a connected rational tree (hence also of simple normal crossing)
and the round up $\lceil L \rceil = L_{\rm red}$.
By the Kawamata-Viehweg vanishing and Riemann-Roch theorem, we have
$h^0(X, K_X + L_{\rm red}) = \frac{1}{2}(K_X + L_{\rm red}) . L_{\rm red} +
\chi(\OO_X) = (-1) + 1 = 0$.
(However, as in the proof of Theorem \ref{theoremk=1} or Lemma \ref{elltype}, we have
$H^0(X, K_X + 2 L_{\rm red}) \ne 0$.)
Therefore, to prove the desired non-vanishing,
one has to show that the above geometric situation will never occur.
\end{remark}

\end{document}